\documentclass[12pt]{article}
\usepackage{amsfonts}
\usepackage{bbm}
\usepackage{mathrsfs}
\usepackage{amssymb}

\def\hang{\hangindent\parindent}
\def\textindent#1{\indent\llap{#1\enspace}\ignorespaces}

\title{On Monoid Graded Local Rings\thanks{Project supported by the National Natural
Science Foundation of China (10971044).}}

\author{Huishi Li\\
{\small Department of Applied Mathematics}\\
{\small College of Information Science and Technology}\\
{\small Hainan University}\\
{\small  Haikou 570228, China}}

\date{}

\begin{document}
\maketitle
\begin{center}
\begin{minipage}{120mm}
{\small {\bf Abstract.} Let $\Gamma$ be a cancelation monoid with
the neutral element $e$. Consider a $\Gamma$-graded ring
$A=\oplus_{\gamma\in\Gamma}A_{\gamma}$, which is not necessarily
commutative. It is proved that $A_e$, the degree-$e$ part of $A$, is
a local ring in the classical sense if and only if the graded
two-sided ideal $\mathfrak{M}$ of $A$ generated by all
non-invertible homogeneous elements is a proper ideal. Defining a
$\Gamma$-graded local ring $A$ in terms of this equivalence, it is
proved that any two minimal homogeneous generating sets of a
finitely generated $\Gamma$-graded $A$-module have the same number
of generators, and furthermore, that most of the basic homological
properties of the local ring $A_e$ hold true for $A$ (at least) in
the $\Gamma$-graded context. }
\end{minipage}\end{center} {\parindent=0pt\par

}\vskip -.5truecm

\renewcommand{\thefootnote}{\fnsymbol{footnote}}
\setcounter{footnote}{-1}


\def\M{\mathfrak{M}}\def\v5{\vskip .5truecm}\def\QED{\hfill{$\Box$}}\def\item{\par\hang\textindent}
\def \r{\rightarrow}\def\OV#1{\overline {#1}}
\def\normalbaselines{\baselineskip 24pt\lineskip 4pt\lineskiplimit 4pt}
\def\mapdown#1{\llap{$\vcenter {\hbox {$\scriptstyle #1$}}$}
                                \Bigg\downarrow}
\def\mapdownr#1{\Bigg\downarrow\rlap{$\vcenter{\hbox
                                    {$\scriptstyle #1$}}$}}
\def\mapright#1#2{\smash{\mathop{\longrightarrow}\limits^{#1}_{#2}}}

\section*{1. Introduction and Preliminaries}
Let $\Gamma$ be a {\it cancelation monoid}, that is, $\Gamma$ is a
semigroup with a neutral element $e$, not necessarily commutative,
and $\Gamma$ satisfies the (left and right) cancelation law:
$\gamma\gamma_1=\gamma\gamma_2$ implies $\gamma_1=\gamma_2$, and
$\gamma_1\gamma =\gamma_2\gamma$ implies $\gamma_1=\gamma_2$ for all
$\gamma ,\gamma_1,\gamma_2\in\Gamma$. This paper studies a
$\Gamma$-graded ring $A=\oplus_{\gamma\in\Gamma}A_{\gamma}$ with the
degree-$e$ part $A_e$ a (commutative or noncommutative) local ring
in the classical sense. As one may see from the literature, such
graded rings may have been widely occurring in both commutative and
noncommutative algebra, commutative and noncommutative algebraic
geometry as well as representation theory, but except for those
well-known classical special cases (see the examples we give after
Definition 2.6 in Section 2), they seemed to have not yet been
explored systematically.  The cut-in point of our work is to prove
that $A_e$ is a local ring in the classical sense if and only if the
graded (two-sided) ideal $\mathfrak{M}$ of $A$ generated by all
non-invertible homogeneous elements is a proper ideal. So it is
natural to call such graded ring $A$ a $\Gamma$-{\it graded local
ring} and ask what properties of the (classical) local ring $A_e$
may hold true for $A$ (at least) at the graded level. Exploring this
question in some detail, although $A$ is not necessarily commutative
and $\M$ is generally not a maximal ideal of the ungraded ring $A$
in the usual sense, and therefore the quotient ring $A/\M$ is
generally not a division ring (or a field) and the localization
techniques as well as the methods dealing with vector spaces no
longer work, we then show that any two minimal homogeneous
generating sets of a finitely generated $\Gamma$-graded $A$-module
have the same number of generators, and furthermore, that most of
the basic homological properties of the local ring $A_e$ hold true
for $A$ (at least) in the $\Gamma$-graded context. Our approach is
carried out through the contents which are arranged as follows.\par

1. Introduction and Preliminaries\par

2. The $\Gamma$-Graded Local Ring $A$\par

3. Minimal Homogeneous Generating Sets of Graded  $A$-Modules\par

4. Finitely Generated gr-Projective $A$-Modules are Free\par

5. The Existence of gr-Projective Covers and Minimal gr-Free
Resolutions over $A$\par

6. Determining Homological Dimensions via $A/\mathfrak{M}$\v5

All rings considered in this paper are associative rings with
multiplicative identity 1, not necessarily commutative; and all
modules over a ring are unitary left modules. By a {\it local ring
in the classical sense} we mean a ring $R$  which has a unique
maximal (left, right and two-sided) ideal ${\bf m}$. It is well
known that if $R$ is a local ring with the maximal ideal {\bf m},
then $r\in R-{\bf m}$ if and only if $r$ is invertible. Commutative
local rings may be seen almost everywhere in each book on
commutative algebra and algebraic geometry. For some interesting
noncommutative local rings we refer to ([Lam], Ch.7, \S 19). \par

Throughout this paper, $\Gamma$ denotes a cancelation monoid with
the neutral element $e$, as we fixed in the beginning. By a
$\Gamma$-graded ring $A$ we mean an associative ring with the
multiplicative identity 1, which has a direct sum decomposition
$A=\oplus_{\gamma\in\Gamma}A_{\gamma}$, where each $A_{\gamma}$ is
an additive subgroup of $A$ (in the case that $A$ is a $K$-algebra
over a field $K$, we require that $A_{\gamma}$ is a $K$-subspace of
$A$), such that $A_{\gamma_1}A_{\gamma_2}\subseteq
A_{\gamma_1\gamma_2}$ for all $\gamma_1,\gamma_2\in\Gamma$. For each
$\gamma\in\Gamma$, we call $A_{\gamma}$ the {\it degree-$\gamma$
part} of $A$. If $0\ne a\in A_{\gamma}$, then we call $a$ a {\it
homogeneous element of degree $\gamma$} and write $d(a)=\gamma$. For
the basics about graded rings and graded modules we used in this
paper, one may refer to [NVO] though it mainly deals with group
$G$-graded rings.\par

Let $A=\oplus_{\gamma\in\Gamma}A_{\gamma}$ be a $\Gamma$-graded ring
as above. By the definition it is clear that the degree-$e$ part
$A_e$ is a subring of $A$, and by using the cancelation law on
$\Gamma$, it is easy to check that $1\in A_e$ (to emphasize its
importance, this property is specified as Lemma 2.1(i) in the next
section). Let $I$ be a (left, right, two-sided) ideal of $A$. If $I$
is generated by homogeneous elements, then we say that $I$ is a {\it
graded} ({\it left, right, two-sided}) {\it ideal} of $A$. By using
the cancelation law on $\Gamma$, it is straightforward to check that
$I$ is a graded (left, right, two-sided) ideal of $A$ if and only if
$I=\oplus_{\gamma\in\Gamma}I_{\gamma}$ with $I_{\gamma}=I\cap
A_{\gamma}$ if and only if $a=\sum_ia_{\gamma_i}\in I$ with
$a_{\gamma_i}\in A_{\gamma_i}$ implies $a_{\gamma_i}\in I$ for all
$i$ if and only if $A/I=\oplus_{\gamma\in\Gamma}(A_{\gamma}+I)/I$.
If each graded left ideal of $A$ is finitely generated, then $A$ is
called a {\it left gr-Noetherian ring}. Similarly, a right
gr-Noetherian ring is defined. If $A$ is both a left and a right
gr-Noetherian ring then it is called a gr-Noetehrian ring.\par

Let $A=\oplus_{\gamma\in\Gamma}A_{\gamma}$ be a $\Gamma$-graded
ring. By a $\Gamma$-graded $A$-module $M$ we mean a left $A$-module
which has a direct sum decomposition
$M=\oplus_{\gamma\in\Gamma}M_{\gamma}$, where each $M_{\gamma}$ is
an additive subgroup of $M$ (in the case that $A$ is a $K$-algebra
over a field $K$, we require that $M_{\gamma}$ is a $K$-subspace of
$M$), such that $A_{\gamma_1}M_{\gamma_2}\subseteq
M_{\gamma_1\gamma_2}$ for all $\gamma_1,\gamma_2\in\Gamma$. For each
$\gamma\in\Gamma$, we call $M_{\gamma}$ the {\it degree-$\gamma$
part} of $M$. If $0\ne \xi\in M_{\gamma}$, then we call $\xi$ a {\it
homogeneous element of degree $\gamma$} and write $d(\xi )=\gamma$.
Let $H$ be a submodule of $M$. If $H$ is generated by homogeneous
elements, then we say that $H$ is a {\it graded submodule} of $M$.
By using the cancelation law on $\Gamma$, it is straightforward to
check that $H$ is a graded submodule of $M$ if and only if
$H=\oplus_{\gamma\in\Gamma}H_{\gamma}$ with $H_{\gamma}=H\cap
M_{\gamma}$ if and only if $\xi =\sum_i\xi_{\gamma_i}\in H$ with
$\xi_{\gamma_i}\in M_{\gamma_i}$ implies $\xi_{\gamma_i}\in H$ for
all $i$ if and only if
$M/H=\oplus_{\gamma\in\Gamma}(M_{\gamma}+H)/H$. If each graded
submodule of $M$ is finitely generated, then $M$ is called a {\it
left gr-Noetherian module}. \par

Let $M=\oplus_{\gamma\in\Gamma}M_{\gamma}$ and
$N=\oplus_{\gamma\in\Gamma}N_{\gamma}$ be $\Gamma$-graded
$A$-modules. If $\varphi$: $M\r N$ is an $A$-module homomorphism
such that $\varphi (M_{\gamma})\subseteq N_{\gamma}$ for all
$\gamma\in\Gamma$, then we call $\varphi$ a {\it graded
$A$-homomorphism}. If $M~\mapright{\varphi}{}~N$ is a graded
$A$-homomorphism, then it is easy to see that Ker$\varphi$ and
Im$\varphi$ are graded submodules of $M$ and $N$ respectively.\par

Note that every $\Gamma$-graded $A$-module $M$ has a generating set
consisting of homogeneous elements, which is usually called a {\it
homogeneous generating set}. The lemma given below tells us that the
cancelation law we assumed on $\Gamma$ guarantees the {\it very
basic graded structure} of  $\Gamma$-graded modules, of which the
verification is straightforward.{\parindent=0pt\v5

{\bf 1.1. Lemma} Let  $A=\oplus_{\gamma\in\Gamma}A_{\gamma}$ be a
$\Gamma$-graded ring and $M=\oplus_{\gamma\in\Gamma}M_{\gamma}$ a
$\Gamma$-graded $A$-module. If $\Omega =\{\xi_i\}_{i\in J}$ is a
homogeneous generating set of $M$, in which $d(\xi_i)=\gamma_i$ (not
necessarily distinct), $i\in J$, then
$$M_{\gamma}=\sum_{i\in J,~\gamma '_i\gamma_i=\gamma}A_{\gamma '_i}\xi_i,~\quad\gamma\in\Gamma .$$\par\QED}\v5

A $\Gamma$-graded $A$-module $F=\oplus_{\gamma\in\Gamma}F_{\gamma}$
is called {\it gr-free} if $F$ is a free $A$-module (in the usual
sense) but with a {\it homogeneous free $A$-basis}, namely $F$ has a
free $A$-basis $\Omega =\{ e_i\}_{i\in J}$ consisting of homogeneous
elements. By Lemma 1.1, if $e_i\in F_{\gamma_i}$ for some
$\gamma_i\in\Gamma$, $i\in J$, then
$$F=\oplus_{i\in J}Ae_i,\quad F_{\gamma}=\sum_{i\in J,~\gamma '_i\gamma_i=\gamma}A_{\gamma '_i}e_i,\quad \gamma\in\Gamma .$$
Indeed, given any free $A$-module $F=\oplus_{i\in J}Ae_i$ and an
arbitrarily chosen subset $\Gamma_1 =\{\gamma_i\}_{i\in J}$ of (not
necessarily distinct) elements of $\Gamma$, a gr-free $A$-module may
be constructed as follows. For each $\gamma\in\Gamma$, if there are
$\gamma_i\in \Gamma_1$, $\gamma '_i\in\Gamma$ such that $\gamma
=\gamma '_i\gamma_i$, then put $F_{\gamma}=\sum_{i\in J,~\gamma_i
'\gamma_i=\gamma}A_{\gamma '_i}e_i$, otherwise put $F_{\gamma}=\{
0\}$. Noticing that $\Gamma$ is a monoid with the neutral element
$e$ and that $\Gamma$ satisfies the cancelation law, one may
directly check that $F=\oplus_{\gamma\in\Gamma}F_{\gamma}$ is the
desired gr-free $A$-module which has the homogeneous free $A$-basis
$\Omega =\{ e_i\}_{i\in J}$ with $e_i\in F_{\gamma_i}$. Thus, if
$M=\oplus_{\gamma\in\Gamma}M_{\gamma}$ is any nonzero
$\Gamma$-graded $A$-module, then, taking a homogeneous generating
set $\{\xi_i~|~d(\xi_i)=\gamma_i\}_{i\in J}$ of $M$, the subset
$\Gamma_1=\{\gamma_i\}_{i\in J}\subseteq\Gamma$ and the gr-free
$A$-module $F=\oplus_{i\in J}Ae_i$ as constructed above by assigning
each $e_i$ the degree $\gamma_i$ , the map $\varphi$: $F\r M$
defined by $\varphi (e_i)=\xi_i$ is a graded $A$-epimorphism. Hence,
every nonzero $\Gamma$-garded $A$-module is a graded homomorphic
image of some gr-free $A$-module.\par

Finally, before introducing the notion of a gr-projective
$A$-module, we also need a fundamental result concerning projective
objects in the $\Gamma$-graded context, which is an analogue of
([NVO], Corollary I.2.2 with $\Gamma =G$ a group; [MR], Proposition
7.6.6 with $\Gamma =\mathbb{N}$ the additive monoid of all
nonnegative integers). As we are working with an arbitrary monoid
$\Gamma$ satisfying the cancelation law, to see why the classical
result holds true as well in our case, a proof is included here.
{\parindent=0pt\v5

{\bf 1.2. Proposition} Let $P=\oplus_{\gamma\in\Gamma}P_{\gamma}$ be
a $\Gamma$-graded $A$-module. The following statements are
equivalent.\par

(i) $P$ is a gr-direct summand of a gr-free $A$-module, i.e.,
$F=P\oplus Q$ in which $F=\oplus_{i\in
J}Ae_i=\oplus_{\gamma\in\Gamma}F_{\gamma}$ is a gr-free $A$-module
and $Q=\oplus_{\gamma\in\Gamma}Q_{\gamma}$ is a graded submodule of
$F$, satisfying $F_{\gamma}=P_{\gamma}+Q_{\gamma}$ for all
$\gamma\in\Gamma$.\par

(ii) Given any graded $A$-epimorphism $\varphi$: $M\r N$ of
$\Gamma$-garded $A$-modules, if $\psi$: $P\r N$ is a graded
$A$-homomorphism, then there is a unique graded $A$-homomorphism
$\OV{\varphi}$: $P\r M$ such that the diagram
$$\begin{array}{ccccc} &&P&&\\
&\scriptstyle{\OV{\varphi}}\swarrow&\mapdownr{\psi}&&\\
M&\mapright{}{\varphi}&N&\r&0\end{array}$$ is commutative.\par

(iii) As an ungraded module, $P$ is a projective $A$-module.\vskip
6pt

{\bf Proof} By virtue of graded $A$-homomorphisms, the proof of (i)
$\Leftrightarrow$ (ii) is similar to the ungraded case. (i)
$\Rightarrow$ (iii) is obvious because, a gr-free $A$-module is
certainly free as an ungraded module.\par

(iii) $\Rightarrow$ (i) Let $\{\xi_i~|~i\in J\}$ be a homogeneous
generating set of $P$, and $\varphi$: $F\r P$ the graded
$A$-epimorphism defined by $\varphi (e_i)=\xi_i$, where
$F=\oplus_{i\in J}Ae_i=\oplus_{\gamma\in\Gamma}F_{\gamma}$ is the
gr-free $A$-module with homogeneous free $A$-basis $\{
e_i~|~d(e_i)=d(\xi_i),~i\in J\}$. Since $P$ is projective as an
ungraded $A$-module, $\varphi$ splits, i.e., there is an
$A$-homomorphism $\beta$: $P\r F$ such that $\varphi\beta =1_P$.
Note that the splitting homomorphism $\beta$ need not be a graded
$A$-homomorphism. But if we define $\OV{\varphi}$: $P\r F$ with
$\OV{\varphi}(x_{\gamma})=f_{\gamma}$, where $x_{\gamma}\in
P_{\gamma}$, $\beta  (x_{\gamma})=f_{\gamma}+\sum_jf_{\gamma_j}$
with $f_{\gamma}\in F_{\gamma}$, $f_{\gamma_j}\in F_{\gamma_j}$ and
$\gamma\ne \gamma_j$, then, with the aid of cancelation law on
$\Gamma$, a direct verification shows that $\OV{\varphi}$ is a
graded $A$-homomorphism, in particular, $\OV{\varphi}(a_{\gamma 
'}x_{\gamma})=a_{\gamma '}\OV{\varphi}(x_{\gamma})$ for all 
$a_{\gamma '}\in A_{\gamma '}$ and $x_{\gamma}\in P_{\gamma}$, 
$\gamma ',\gamma\in\Gamma$ (this is the {\it key} point), such that 
$\varphi\OV{\varphi}=1_P$. Thus, the graded $A$-homomorphism 
$\OV{\varphi}$ splits $\varphi$ in degrees and therefore, (i) is 
proved.} \QED \v5

A $\Gamma$-graded $A$-module $P=\oplus_{\gamma\in\Gamma}P_{\gamma}$
is called a {\it gr-projective module} if it satisfies one of the
equivalent conditions of Proposition 1.2.

\section*{2. The $\Gamma$-Graded Local Ring $A$}
Let $\Gamma$ be a cancelation monoid with the neutral element $e$,
and  $A=\oplus_{\gamma\in\Gamma}A_{\gamma}$ a $\Gamma$-graded ring.
In this section we introduce the notion of a $\Gamma$-graded local
ring by showing that the degree-$e$ part $A_e$ of $A$ is a local
ring in the classical sense if and only if the graded two-sided
ideal $\M$ of $A$ generated by all non-invertible homogeneous
elements is a proper ideal.\v5

We start with some preliminaries. First note that $A_e$ is a subring
of $A$. Let $a\in A_{\gamma}$ be a homogeneous element of degree
$\gamma$. We say that $a$ is {\it left} ({\it right}) {\it
gr-invertible} if there is a homogeneous element $b$ of $A$ such
that $ba=1$ ($ab=1$). If $a$ is both left and right gr-invertible,
we say that $a$ is {\it gr-invertible}. The properties (i) and (ii)
mentioned in the first lemma below may be directly checked by using
the cancelation law assumed on $\Gamma$. {\parindent=0pt\v5

{\bf 2.1. Lemma}  (i) The (multiplicative) identity 1 of $A$ is a
homogeneous element of degree $e$, i.e., $1\in A_e$.\par

(ii) If a homogeneous element $a\in A_{\gamma}$ is left (right)
invertible in $A$, then it is left (right) gr-invertible.\par
\QED}\v5

From now on in our discussion we shall freely use the property Lemma
2.1(i) without additional indication. \v5

Now, suppose that the degree-$e$ part $A_e$ of $A$ is a local ring
in the classical sense. Writing {\bf m} for the unique maximal ideal
of $A_e$, as usual we use $(A_e,{\bf m})$ to indicate that $A_e$ is
a local ring with the maximal ideal {\bf m}. {\parindent=0pt\v5

{\bf 2.2. Lemma} Let $A=\oplus_{\gamma\in\Gamma}A_{\gamma}$ and
$(A_e,{\bf m})$ be as fixed above. If a homogeneous element
$a_{\gamma}\in A_{\gamma}$ is left invertible, then it is also right
invertible and hence invertible. Similarly, if a homogeneous element
$a_{\gamma}\in A_{\gamma}$ is right invertible, then it is also left
invertible and hence invertible.\vskip 6pt

{\bf Proof} Suppose that the homogeneous element $a_{\gamma}\in
A_{\gamma}$ is left invertible. Then, it is left gr-invertible by
Lemma 2.1, i.e., there is a homogeneous element, say
$b_{\gamma_1}\in A_{\gamma_1}$, such that
$b_{\gamma_1}a_{\gamma}=1$. Let $a_{\gamma}b_{\gamma_1}=c$. Then $c$
is a homogeneous element and
$$a_{\gamma}=a_{\gamma}(b_{\gamma_1}a_{\gamma})=ca_{\gamma}.\eqno{(1)}$$
Since $\Gamma$ satisfies the right cancelation law, it follows from
the equality (1) that $c\in A_e$. Note that from (1) we also have
$(1-c)a_{\gamma}=0$. Noticing  $a_{\gamma}\ne 0$, we conclude that
$c\not\in {\bf m}$. Thus, $c$ is invertible in the local ring $A_e$
and hence $a_{\gamma}(b_{\gamma_1}c^{-1})=1$, showing that
$a_{\gamma}$ is also right invertible.}\par

By using the left cancelation law on $\Gamma$, the second assertion
is proved in the same way.\par\QED{\parindent=0pt \v5

{\bf 2.3. Proposition} Let $A=\oplus_{\gamma\in\Gamma}A_{\gamma}$
and $(A_e,{\bf m})$ be as in Lemma 2.2. The following statements
hold.\par

(i) Let $\{ L_i~|~i\in I\}$ be the set of all proper graded left
ideals of $A$ and $\mathfrak{L}=\sum_{i\in I} L_i$.  Then
$\mathfrak{L}$ is a proper graded left ideal of $A$.\par

(ii) Let $\{ R_j~|~j\in J\}$ be the set of all proper graded right
ideals of $A$ and $\mathfrak{R}=\sum_{j\in J} R_j$. Then
$\mathfrak{R}$ is a proper graded right ideal of $A$.\par

(iii) Let $\{ T_k~|~k\in K\}$ be the set of all proper graded
two-sided ideals of $A$ and $\mathfrak{M}=\sum_{k\in K}T_k$. Then
$\mathfrak{M}$ is a proper graded two-sided ideal of $A$.\par

(iv) $\mathfrak{L}=\mathfrak{M}=\mathfrak{R}$. \par

(v) Each homogeneous element $a\in A_{\gamma}-\mathfrak{M}$ is
gr-invertible, $\gamma\in\Gamma$, namely there is a homogeneous
element $b$ of $A$ such that $ab=ba=1$. Hence, $\M$ is, as a graded
ideal of $A$, generated by all non-invertible homogeneous elements
of $A$.\par

(vi)  $\mathfrak{M}\cap A_e={\bf m}$.\vskip 6pt

{\bf Proof} (i) -- (iii) Suppose the contrary that $1\in
\mathfrak{L}$. Then $1=\sum^m_{t=1}a_{i_t}$ with $a_{i_t}\in
(L_{i_t})_e=L_{i_t}\cap A_e$, and there is at least one $a_{i_t}$
which is not contained in the maximal ideal {\bf m} of $A_e$. But
then $a_{i_t}$ is invertible in $A_e$ and thus $L_{i_t}$ would not
be a proper graded left ideal, a contradiction. Hence,
$\mathfrak{L}$ is a proper graded left ideal of $A$. In a similar
way, we can prove that $\mathfrak{R}$ is a proper graded right ideal
of $A$, and that $\M$ is a proper graded two-sided ideal of $A$.\par

(iv) We need only to show that $\mathfrak{L}=\mathfrak{R}$ because,
by (i) -- (iii),  it is then clear that
$\mathfrak{L}=\mathfrak{M}=\mathfrak{R}$. Let $a\in
\mathfrak{L}_{\gamma}$ be a nonzero homogeneous element of degree
$\gamma$ in $\mathfrak{L}$. We claim that the graded right ideal
$aA$ is a proper right ideal of $A$ and hence $a\in aA\subset
\mathfrak{R}$. Otherwise, $1\in aA$ would mean that $a$ is right
invertible and then, it follows from Lemma 2.2 that $a$ is also left
invertible, implying that $\mathfrak{L}$ would not be a proper left
ideal. This contradiction confirms our claim. Since  $a\in aA\subset
\mathfrak{R}$ holds for each homogeneous element $a$ of
$\mathfrak{L}$, we have $\mathfrak{L}\subseteq \mathfrak{R}$.
Similarly, by using  Lemma 2.2 we can prove $\mathfrak{R}\subseteq
\mathfrak{L}$. Therefore, $\mathfrak{L}=\mathfrak{R}$, as
desired.\par

(v) If $a\in A_{\gamma}-\mathfrak{M}$, then it follows from (iv)
that $Aa=A=aA$. Hence $a$ is gr-invertible by Lemma 2.1.\par

(vi) Since $\mathfrak{M}$ is a proper graded two-sided ideal of $A$
by (iii), each $ a\in\mathfrak{M}\cap A_e$ is not invertible. So,
$a\in {\bf m}$. This shows that $\mathfrak{M}\cap A_e\subseteq {\bf
m}$. Conversely, for each $0\ne b\in {\bf m}$, we claim that $Ab$ is
a proper graded left ideal of $A$ and hence $b\in
Ab\subset\mathfrak{M}$ by (iv). Otherwise, $1\in Ab$ would mean that
$b$ is invertible in $A_e$, contradicting $b\in {\bf m}$. Thus, our
claim is true and hence ${\bf m}\subseteq\M\cap A_e$. Consequently,
$\mathfrak{M}\cap A_e={\bf m}$.\QED}\v5

Conversely, we have the next proposition, of which the verification
is straightforward and so is omitted. {\parindent=0pt\v5

{\bf 2.4. Proposition} Let $A=\oplus_{\gamma\in\Gamma}A_{\gamma}$ be
a $\Gamma$-graded ring, and let $\M$ be the graded two-sided ideal
of $A$ generated by all non-invertible homogeneous elements. If $\M$
is a proper ideal of $A$, then the following statements hold.\par

(i) Each homogeneous element $a\in A_{\gamma}-\M$ is gr-invertible,
$\gamma\in \Gamma$. Hence, $\M$ is the unique maximal graded (left,
right and two-sided) ideal of $A$. \par

(ii) $A_e$ is a local ring in the classical sense, which has the
unique maximal (left, right and two-sided) ideal ${\bf m}=\M \cap
A_e$.\par

(iii) The quotient ring $D=A/\M$ is a $\Gamma$-graded division ring,
that is, each nonzero homogeneous element of $D$ is gr-invertible.
In particular, $D_e=A_e/{\bf m}$ is a division ring.\par\QED}\v5

Combining Proposition 2.3 and Proposition 2.4, we are able to
mention the main result of this section. {\parindent=0pt\v5

{\bf 2.5. Theorem} Let $\Gamma$ be a cancelation monoid with the
neutral element $e$, and let $A=\oplus_{\gamma\in\Gamma}A_{\gamma}$
be a $\Gamma$-graded ring. Then $A_e$, the degree-$e$ part of $A$,
is a local ring in the classical sense if and only if the graded
two-sided ideal $\M$ of $A$ generated by all non-invertible
homogeneous elements is a proper ideal. \QED}\v5

Naturally, Theorem 2.5 leads to the following{\parindent=0pt\v5

{\bf 2.6. Definition} Let $\Gamma$ be a cancelation monoid with the
neutral element $e$, and let $A=\oplus_{\gamma\in\Gamma}A_{\gamma}$
be a $\Gamma$-graded ring. If $A_e$, the degree-$e$ part of $A$, is
a local ring in the classical sense, then we call $A$ a
$\Gamma$-{\it graded local ring}. }\v5

By Definition 2.6, the class of $\Gamma$-graded local rings includes
a lot of well-known algebras in both the commutative and
noncommutative contexts. For instance, any (commutative or
noncommutative) connected $\mathbb{N}$-graded $K$-algebra
$A=\oplus_{p\ge 0}A_p$ with $A_0=K$ a field; any (commutative or
noncommutative) $\mathbb{N}$-graded ring $A=\oplus_{p\ge 0}A_p$ with
$A_0$ a local ring (in the commutative case such $A$ is called a
{\it generalized local ring} in the sense of [GW], see also [Eis]
P.510), as to which, some results concerning the global homological 
 dimension were obtained in [Li1]; and any commutative 
$\mathbb{Z}$-graded (or $\mathbb{N}$-graded)
*local ring in the sense of [BH].\par

More generally, let $\Gamma$ be a cancelation monoid with the
neutral element $e$, and $A=\oplus_{\gamma\in\Gamma}A_{\gamma}$ a
$\Gamma$-graded ring such that $A_e$ is a commutative ring contained
in the center of $A$. If $p$ is a prime ideal of $A_e$, then, with
respect to the multiplicative subset $S_p=A_e-p$,  the classical
central localization
$S_p^{-1}A=\oplus_{\gamma\in\Gamma}S_p^{-1}A_{\gamma}$ of $A$ at $p$
is a $\Gamma$-graded local ring in which the degree-$e$ part
$(S_p^{-1}A)_e=S_p^{-1}A_e$ is a local ring with the maximal ideal
${\bf m}=S^{-1}_pp$; if furthermore $A_e$ is a commutative domain
and $S=A_e-\{ 0\}$, then $(S^{-1}A)_e=S^{-1}A_e$ is a field.
{\parindent=0pt\v5

{\bf Convention} In light of Proposition 2.3 and Proposition 2.4,
from now on in this paper we shall also alternatively use the
equivalent of Definition 2.6, that is, we say that the
$\Gamma$-graded ring $A=\oplus_{\gamma\in\Gamma}A_{\gamma}$ is a
$\Gamma$-graded local ring if the graded two-sided ideal $\M$ of $A$
generated by all non-invertible homogeneous elements is the unique
maximal graded (left, right and two-sided) ideal of $A$, and we just
simply say that $\M$ is the maximal graded ideal of $A$. Moreover,
the property ${\bf m}=\M\cap A_e$ will be freely used without extra
indication, where {\bf m} is the unique maximal (left, right and
two-sided) ideal of $A_e$.\v5

{\bf Remark} Let $=\oplus_{\gamma\in\Gamma}A_{\gamma}$ be a
$\Gamma$-garded local ring with the maximal graded ideal $\M$. Then
$D=A/\M$ is a $\Gamma$-graded division ring by Proposition 2.4. If
$\Gamma =\mathbb{Z}$ is the additive group of all integers, then it
follows from ([NVO], Corollary I.4.3) that either $D=D_0$ in case
the gradation is trivial, or else $D\cong D_0[x,x^{-1},\varphi]$,
where $\varphi$ is a ring automorphism of $D_0$, and $x$ is an
indeterminate such that $x\lambda =\varphi (a)\lambda$ for all
$\lambda\in D_0$; moreover, $D$ is a left and right principal ideal
domain. However, as pointed out in ([NVO], Remark I.4.4), for more
general monoid $\Gamma$ (even if $\Gamma$ is a group) the graded
division ring $D$ may not necessarily have the structure as with
$\Gamma =\mathbb{Z}$. Nevertheless, some basic properties of modules
over $D$ will be discussed in Section 4. }

\section*{3. Minimal Homogeneous Generating Sets of Graded $A$-Modules}\par

Throughout this section, we let $\Gamma$ be a cancelation monoid
with the neutral element $e$, and let
$A=\oplus_{\gamma\in\Gamma}A_{\gamma}$ be a $\Gamma$-graded local
ring (in the sense of Definition 2.6) with the maximal graded ideal
$\M$. Our aim below is to show that any two minimal homogeneous
generating sets of a finitely generated $\Gamma$-graded $A$-module
$M$ have the same number of generators, and that in certain cases,
the number of degree-$\gamma$ homogeneous elements contained in each
minimal homogeneous generating set is the same.\v5

Let $M$ be a $\Gamma$-graded $A$-module with homogeneous generating
set $\Omega =\{\xi_i\}_{i\in J}$. If any proper subset of $\Omega$
cannot be a generating set of $M$, then we call $\Omega$ a {\it
minimal homogeneous generating set} of $M$.  It is easy to see that
any homogeneous free $A$-basis $X$ of a gr-free $A$-module $F$ is a
minimal homogeneous generating set of $F$. But for an arbitrary
$A$-module $M$, if $M$ is not finitely generated, in general we do
not know whether $M$ has a (finite or infinite) minimal generating
set or not. Nevertheless, for the completeness, let us first make it
clear that if an $A$-module has an infinite minimal generating set,
then any two minimal generating sets of $A$ have the same
cardinality. The argument on this property is actually the same as
dealing with vector spaces over a (skew) field (see e.g. [Lan]). To
be convinced, we also present a detailed proof.
\par

We first mention an easy but useful lemma.

{\parindent=0pt\v5

{\bf 3.1. Lemma} Let $M$ be an $A$-module. If $M$ is finitely
generated, then any minimal generating set of $M$ is finite.
\par\QED\v5

{\bf 3.2. Theorem}  Let $M$ be a left $A$-module  with an infinite
minimal generating set $T$. Then any other minimal generating set
$T'$ of $M$ has the same cardinality as $T$.\vskip 6pt

{\bf Proof} First note that $T'$ must be infinite by Lemma 3.1.
Since $T$ generates $M$, for each $\eta\in T'$, we have $\eta$
contained in the submodule generated by a finite subset $T_{\eta}$
of $T$. Put $T^*=\cup_{\eta\in T'}T_{\eta}$. Then $T'\subset\langle
T^*\rangle$ and hence $\langle T'\rangle =M=\langle T^*\rangle$.
Thus $T^*=T$ by the minimality of $T$. Furthermore, as each
$T_{\eta}$ is finite and $T'$ is infinite, we can find an injective
map from $T_{\eta}$ into $T'$. It follows that we can find an
injective map from $T^*$ to $T'\times T'$. So, the cardinality $|T|$
of $T=T^*$ is at most the cardinality $|T'\times T'|$ of $T'\times
T'$. But since $T'$ is infinite, $|T'|=|T'\times T'|$ holds by
Theorem 3.6 from Appendix 2 of [Lan]. Consequently $|T|\le |T'|$. By
symmetry we conclude $|T|=|T'|$.\QED}\v5

We now deal with finitely generated $\Gamma$-graded $A$-modules.
{\parindent=0pt\v5

{\bf 3.3. Proposition} Let $M$ be a $\Gamma$-graded $A$-module and
let $T_1=\{ \xi_1,\ldots ,\xi_m, \beta_1,\ldots ,\beta_q\}$ and
$T_2=\{ \eta_1,\ldots ,\eta_n, \beta_1,\ldots ,\beta_q\}$ be two
finite homogeneous generating sets of $M$ such that $\{ \xi_1,\ldots
,\xi_m\}\cap\{ \eta_1,\ldots ,\eta_n\}=\emptyset$. Suppose that
$T_2$ is a minimal homogeneous generating set of $M$. The following
two statements hold.\par

(i) $n\le m$.\par

(ii) Let $n_i$ denote the number of homogeneous elements of degree
$\gamma_i$ in $T_2$. If $\Gamma$ is a totally ordered monoid with
the total ordering $\prec$, such that the neutral element $e$ of
$\Gamma$ is the smallest element of $\Gamma$, then $T_1$ contains at
least $n_i$ homogeneous elements of degree $\gamma_i$.\par \vskip
6pt

{\bf Proof} (i) Suppose the contrary that $n>m$. Since both $T_1$,
$T_2$ are homogeneous generating sets of $M$,  the $\xi_i$'s can be
generated by $T_2$, that is, there are homogeneous elements $B_{ij},
 D_{ir}\in A$ such that
$$\begin{array}{rcl} \xi_i&=&\sum^n_{r=1}D_{ir}\eta_r+\sum^q_{j=1}B_{ij}\beta_j\\
&{~}&\hbox{with}~d(\xi_i)=d(D_{ir}\eta_r)=d(B_{ij}\beta_j);\end{array}1\le
i\le m\eqno{(1)}$$ and on the other hand, there are homogeneous
elements $H_{1i},E_{1j}\in A$ such that
$$\begin{array}{rcl} \eta_1&=&\sum^m_{i=1}H_{1i}\xi_i+\sum^q_{j=1}E_{1j}\beta_j\\
&{~}&\hbox{with}~d(\eta_1)=d(H_{1i}\xi_i)=d(E_{1j}\beta_j).\end{array}\eqno{(2)}$$
Substituting each $\xi_i$ in (2) by (1), we have
$$\begin{array}{rcl} \eta_1&=&\displaystyle{\sum^m_{i=1}}H_{1i}\left (\sum^q_{j=1}B_{ij}\beta_j+
\displaystyle{\sum^n_{r=1}}D_{ir}\eta_r\right )
+\displaystyle{\sum^q_{j=1}}E_{1j}\beta_j\\
\\
&=&\displaystyle{\sum^q_{j=1}}\left (\left
(\displaystyle{\sum^m_{i=1}}H_{1i}B_{ij}\beta_j\right
)+E_{1j}\beta_j\right
)+\displaystyle{\sum^n_{r=2}\sum_{i=1}^m}H_{1i}D_{ir}\eta_r+
\displaystyle{\sum^m_{i=1}}H_{1i}D_{i1}\eta_1.\end{array}\eqno{(3)}$$
Put $f=\sum^m_{i=1}H_{1i}D_{i1}$. Then $f\ne 0$. Otherwise, (3)
would mean that $T_2$ is not minimal. Since all homogeneous terms on
the right-hand side of (3) have the same degree $d(\eta_1)$ by (1)
and (2), the right cancelation law we assumed on $\Gamma$ entails
that all homogeneous elements $H_{1i}D_{i1}$ occurring in the
representation of $f$  must be contained in $A_e$. Note that $A_e$
is a local ring in the classical sense, which has the maximal ideal
${\bf m}=\M\cap A_e$ (Proposition 2.4). So, $1-v$ is invertible in
$A_e$ for all $v\in {\bf m}$. It follows once again from the
minimality of $T_2$ that there is at least one $i^*$ such that
$H_{1i^*}D_{i^*1}\ne 0$ and $H_{1i^*}D_{i^*1}\in A_e-{\bf m}$.
Consequently, there is some $c\in A_e-{\bf m}$ such that
$H_{1i^*}(D_{i^*1}c)=1.$ This shows that $H_{1i^*}$ is right
invertible. But by Lemma 2.2, $H_{1i^*}$ is also left invertible and
hence invertible. Therefore, rewriting (2) as
$$\eta_1=H_{1i^*}\xi_{i^*}+\sum^m_{i\ne i^*}H_{1i}\xi_i+\sum^q_{j=1}E_{1j}\beta_j,$$
we have
$$\xi_{i^*}=H_{1i^*}\mathbf{}^{-1}\left (\eta_1-\sum^m_{i\ne i^*}H_{1i}\xi_i-
\sum^q_{j=1}E_{1j}\beta_j\right ) .\eqno{(4)}$$ Since $T=T_1\cup\{
\eta_1\}$ is also a homogeneous generating set for $M$,  the formula
(4) turns out that
$$T_3=T-\{ \xi_{i^*}\}=\{ \xi_1,\ldots ,\xi_{i^*-1},\xi_{i^*+1},\ldots ,\xi_m,
\beta_1,...,\beta_q,\eta_1\}$$ is a homogeneous generating set of
$M$.}\par

With $T_1$ replaced by $T_3$, we can repeat the substitution
procedure as presented above and so forth, until the remained
$\xi_1,\ldots ,\xi_{i^*-1},\xi_{i^*+1},\ldots $ and $\xi_m$ are
replaced by $\eta_2,\eta_3,\ldots $ and $\eta_{m}$ successively. So,
totally after $m$ substitutions we reach a homogeneous generating
set $\{ \eta_1,\ldots ,\eta_{m},\beta_1,\ldots ,\beta_q\}$ for $M$,
which is a proper subset of $T_2$. This contradicts the minimality
of $T_2$. Therefore $n\le m$, as desired.{\parindent=0pt\par

(ii) If $\Gamma$ is a totally ordered monoid with the total ordering
$\prec$, then  $\gamma_1\prec \gamma_2$ implies $\gamma_1\gamma\prec
\gamma_2\gamma$ and $\gamma\gamma_1\prec\gamma\gamma_2$ for all
$\gamma_1,\gamma_2,\gamma\in\Gamma$. By the totality of $\prec$ it
is easy to see that $\Gamma$ satisfies the (left and right)
cancelation law. If furthermore the neutral element $e$ of $\Gamma$
is the smallest element, then it is clear that the homogeneous
element $H_{1i^*}$, obtained in the proof of (i) such that
$H_{1i^*}D_{i^*1}\in A_e-{\bf m}$, itself is contained in $A_e-{\bf
m}$. It follows that the inverse element $H_{1i^*}^{-1}$ of 
$H_{1i^*}$ is  of degree $e$ and hence $d(\xi_{i^*})=d(\eta_1)$. 
This proves the desired assertion.\QED }\v5

We are ready to mention the main result of this
section.{\parindent=0pt\v5

{\bf 3.4. Theorem} Let $A$ be a $\Gamma$-graded local ring as before
and let $M$ be a finitely generated $\Gamma$-graded $A$-module.\par

(i) Any two minimal homogeneous generating sets of $M$ have the same
number of generators.\par

(ii) If $\Gamma$ is a totally ordered monoid with the total ordering
$\prec$, such that the neutral element $e$ of $\Gamma$ is the
smallest element of $\Gamma$, then any two minimal homogeneous
generating sets of $M$ contain the same number of homogeneous
elements of degree $\gamma$ for all $\gamma\in\Gamma$.\QED}\v5

We end this section by pointing out that in a similar way as
presented in the proof of Proposition 3.3, it may be shown that if
$A_e$ is contained in the center of $A$, then the results of Theorem
3.4 hold true for any finitely generated graded two-sided ideal of
$A$, that is, we can mention the next{\parindent=0pt\v5

{\bf 3.5. Proposition} Let $A$ be a $\Gamma$-graded local ring as
before and let $I$ be a finitely generated graded two-sided ideal of
$A$. Suppose that $A_e$ is contained in the center of $A$, then the
following two statements hold.\par

(i) Any two minimal homogeneous generating sets of $I$ have the same
number of generators.\par

(ii) If $\Gamma$ is a totally ordered monoid with the total ordering
$\prec$, such that the neutral element $e$ of $\Gamma$ is the
smallest element of $\Gamma$, then any two minimal homogeneous
generating sets of $I$ contain the same number of homogeneous
elements of degree $\gamma$ for all $\gamma\in\Gamma$.}\v5

\section*{4. Finitely Generated gr-Projective $A$-Modules are Free}
Let $\Gamma$ be a cancelation monoid with the neutral element $e$,
and let $A=\oplus_{\gamma\in\Gamma}A_{\gamma}$ be a $\Gamma$-graded
local ring (in the sense of Definition 2.6) with the maximal graded
ideal $\M$. In this section we show that every finitely generated
gr-projective $A$-module is free. Concerning the basics of gr-free
$A$-modules and gr-projective modules, one is referred to Section
1.\v5

Noticing that every $\Gamma$-graded division ring is a special
$\Gamma$-graded local ring, we start by showing that $\Gamma$-graded
modules over a $\Gamma$-graded division ring
$D=\oplus_{\gamma\in\Gamma}D_{\gamma}$ behave very like modules over
the usual (ungraded) division rings. Since each nonzero homogeneous
element of $D$ is gr-invertible, i.e., if $0\ne s_{\gamma}\in
D_{\gamma}$ then there is some homogeneous element, say
$s_{\gamma_1}\in D_{\gamma_1}$, such that
$s_{\gamma_1}s_{\gamma}=1=s_{\gamma}s_{\gamma_1}$, given a nonzero
$\Gamma$-graded $D$-module $E$, naturally we can talk about the
gr-linear independence of homogeneous elements in $E$. More
precisely, let $T$ be a nonempty subset of $E$ consisting of
homogeneous elements. For any $q\ge 1$, $\xi_{i_1},\ldots
,\xi_{i_q}\in T$  and homogeneous elements $s_{\gamma_1},\ldots
,s_{\gamma_q}\in D$, if $d(s_{\gamma_1}\xi_{i_1})=\cdots
=d(s_{\gamma_q}\xi_{i_q})$ and
$\sum_{j=1}^qs_{\gamma_j}\xi_{i_j}=0$ imply $s_{\gamma_1}=\cdots
=s_{\gamma_q}=0$, then  we say that $T$ is a {\it gr-linearly
independent homogeneous subset} of $E$. If, with respect to the set
theoretic inclusion relation, $T$ is maximal among all gr-linearly
independent homogeneous subsets in $E$,  then we call $T$ a {\it
maximal gr-linearly independent homogeneous subset}.
\par

By the cancelation law on $\Gamma$ it is easy to see that the
gr-linear independence of homogeneous elements defined above is
equivalent to the usual linear independence of homogeneous elements
over $D$, namely for any $q\ge 1$, $\xi_{i_1},\ldots ,\xi_{i_q}\in
T$  and $f_1,\ldots ,f_q\in D$, if $\sum_{j=1}^qf_j\xi_{i_j}=0$ then
$f_1=\cdots =f_t=0$. {\parindent=0pt\v5

{\bf 4.1. Proposition} Let $D=\oplus_{\gamma\in\Gamma}D_{\gamma}$ be
a $\Gamma$-graded division ring, and $E$ a nonzero $\Gamma$-graded
$D$-module. The following statements hold.\par

(i) Each gr-linearly independent homogeneous subset $V$ of $E$ is
contained in a maximal gr-linearly independent homogeneous subset of
$E$. \par

(ii) $E$ is a gr-free $D$-module (and hence a free $D$-module).
Indeed, each maximal gr-linearly independent homogeneous subset $S$
of $E$ forms a homogeneous free $D$-basis for $E$;  and moreover,
$S$ also plays the role as a minimal homogeneous generating set of
$E$.\par

(iii) Each  minimal homogeneous generating set $S$ of $E$ is a
homogeneous free $D$-basis for $E$, and moreover, $S$ is also a
maximal gr-linearly independent homogeneous subset of $E$. \par

(iv) Each homogeneous free $D$-basis of $E$ is a minimal homogeneous
generating set of $E$. Any two homogeneous free $D$-bases of $E$
have the same cardinality. \vskip 6pt

{\bf Proof} (i) Let $V$ be a gr-linearly independent homogeneous
subset of $E$ and $\Sigma$ the set of all gr-linearly independent
homogeneous subsets containing $V$ in $E$. Since $E\ne \{ 0\}$ and
each nonzero homogeneous element in $D$ is gr-invertible, we have
$\Sigma\ne\emptyset$. Partially order $\Sigma$ by set theoretic
inclusion. Then, actually as for establishing the existence of a
basis for vector spaces over a division ring, we can use Zorn's
lemma to find a maximal element $S$ in $\Sigma$ which is clearly a
maximal  gr-linearly independent homogeneous subset of $E$.\par

(ii) It is straightforward to check that a maximal gr-linearly
independent homogeneous subset $S$ is a homogeneous free $D$-basis
of $E$, and that $S$ is also a minimal homogeneous generating set of
$E$.\par

(iii) The existence of a minimal homogeneous generating set for $E$
is given by (i) and (ii). Let $S$ be any minimal homogeneous
generating set of $E$. Then it is straightforward to verify that $S$
is a gr-linearly independent homogeneous subset of $E$ and hence a
homogeneous free $D$-basis for $E$. To see $S$ is a maximal
gr-linearly independent homogeneous subset, let $S_1$ be a
gr-linearly independent homogeneous subset such that $S\subseteq
S_1$. But by (i) and (ii), there is a maximal gr-linearly
independent homogeneous subset $S_2$, which is also a minimal
homogeneous generating set of $E$,  such that $S\subseteq
S_1\subseteq S_2$. It follows that $S=S_1=S_2$, showing that $S$ is
maximal. \par

(iv) That a homogeneous free $D$-basis of $E$ is a minimal
homogeneous generating set of $E$ is clear. Since a $\Gamma$-graded
division ring is certainly a $\Gamma$-graded local ring in the sense
of Definition 2.6, from the discussion of Section 3 we know that any
two minimal homogeneous generating sets of $E$ have the same
cardinality. \QED }\v5

Let $D=\oplus_{\gamma\in\Gamma}D_{\gamma}$ and $E$ be as in
Proposition 4.1. Then we can define the {\it gr-rank} of $E$ as the
cardinality of a minimal homogeneous generating set of $E$.\v5

Now, let us return to a general $\Gamma$-graded local ring
$A=\oplus_{\gamma\in\Gamma}A_{\gamma}$ with the maximal graded ideal
$\M$. Then, the degree-$e$ part $A_e$ of $A$, where $e$ is the
neutral element of $\Gamma$, is a local ring (in the classical
sense) with the maximal ideal ${\bf m}=\M \cap A_e$. Moreover, we
note that $\M$ is the graded Jacobson radical $J^g(A)$ of $A$ in the
usual graded context (c.f. [NVO]), i.e., $J^g(A)=\mathfrak{M}$. So,
a monoid graded version of Nakayama's lemma holds true for
$\Gamma$-graded $A$-modules. To see how $\M$ works in the monoid
graded case, we include a proof here. {\parindent=0pt\v5

{\bf 4.2. Lemma} Let $L$ be a graded left ideal of $A$ which is
contained in $\M$, and let $M$ be a finitely generated
$\Gamma$-graded $A$-module. If $LM=M$, then $M=\{ 0\}$. \vskip 6pt

{\bf Proof} We do induction on the number of generators of $M$. Let
$\xi_1,\ldots ,\xi_i$ be homogeneous generators of $M$, i.e.,
$M=\sum^s_{i=1}A\xi_i$. If $LM=M$, then by the cancelation law on
$\Gamma$ (or Lemma 1.1), $\xi_1$ has a representation
$$\xi_1=c_1\xi_1+c_2\xi_2+\cdots +c_s\xi_s,$$
in whcih each $c_i$ is a homogeneous element of $LA\subseteq\M$ such
that $d(c_i\xi_i)=d(\xi_1)$ for each $c_i\xi_i\ne 0$. If $c_1=0$,
then $M$ can be generated by $s-1$ elements. If $c_1\ne 0$, then by
the cancelation law on $\Gamma$ we have $c_1\in\M\cap A_e={\bf m}$.
It follows that  $\xi_1=(1-c_1)^{-1}\sum^s_{i= 2}c_i\xi_i$, showing
that $M$ can be generated by $s-1$ elements. Hence the proof is
complete by induction.\QED\v5

{\bf 4.3. Corollary} Let $L$ be a graded left ideal of $A$ which is
contained in $\M$, and let $M$ be a finitely generated
$\Gamma$-graded $A$-module, and $H$ a graded submodule of $M$. If
$M=H+LM$, then $M=H$.\par\QED\v5

{\bf 4.4. Proposition} Let $M$ be a nonzero and finitely generated
$\Gamma$-graded $A$-module, and $\Omega =\{\xi_1,\ldots
,\xi_s\}\subset M$ a subset of homogeneous elements. Considering the
$\Gamma$-graded division ring $D=A/\M$, the following two statements
are equivalent.\par

(i) $\Omega$ is a minimal homogeneous generating set of $M$.\par

(ii) $\OV{\Omega}=\{\OV{\xi_1},\ldots ,\OV{\xi_s}\}$, where each
$\OV{\xi_i}=\xi_i+\M M$ is the residue class of $\xi_i$ in $ E=M/\M
M$, is a homogeneous $D$-basis for the free $D$-module $E$.\vskip
6pt

{\bf Proof} (i) $\Rightarrow$ (ii) Since $M\ne \{ 0\}$ and hence
$E\ne \{ 0\}$ by Lemma 4.2, that $\OV{\Omega}$ generates $E$ is 
clear. It remains to show that $\OV{\Omega}$ is a gr-linearly 
independent homogeneous subset in $E$. To this end, by Lemma 1.1 we 
assume that $\sum_{i=1}^s\OV{a_i}\OV{\xi_i}=0$, in which each 
$\OV{a_i}=a_i+\M\in D=A/\M$ with $a_i\in A-\M$ a homogeneous 
element, and all the  $\OV{a_i}\OV{\xi_i}$ have the same degree, say 
$\gamma$. Then $\sum^s_{i=1}a_i\xi_i\in\M M\cap M_{\gamma}$ and 
$\sum^s_{i=1}a_i\xi_i=\sum^s_{i=1}b_i\xi_i$, where, by the 
cancelation law on $\Gamma$, each $b_i$ is a homogeneous element in 
$\M$ and all the nonzero $b_i\xi_i$ have the same degree $\gamma$. 
Noticing that $A$ is $\Gamma$-graded local and hence each 
homogeneous $a_i\in A-\M$ is gr-invertible, we should have 
$a_i\xi_i\ne0$, $1\le i\le s$. Since $a_1\xi_1\ne 0$ implies
$$\xi_1=a_1^{-1}b_1\xi_1+a_1^{-1}\sum^s_{i=2}(b_i-a_i)\xi_i,$$
the minimality of $\Omega$ entails  $a_1^{-1}b_1\xi_1\ne 0$. Thus,
the homogeneous element $\xi_1$ and the homogeneous element
$a_1^{-1}b_1\xi_1$ have the same degree.  It follows from the
cancelation law on $\Gamma$ that $a_1^{-1}b_1\in A_e\cap\M ={\bf
m}$, and consequently,
$\xi_1=(1-a_1^{-1}b_1)^{-1}a_1^{-1}\sum^s_{i=2}(b_i-a_i)\xi_i,$
contradicting the minimality of $\Omega$. This shows that we must
have $a_1\in\M$, i.e., $\OV{a_1}=0$. Similarly, $\OV{a_i}=0$ for
$2\le i\le s$. Therefore $\OV{\Omega}$ is a gr-linearly independent
homogeneous subset of $E$, as desired. \par

(ii) $\Rightarrow$ (i) Suppose that $\OV{\Omega}$ is a $D$-basis for
the free $D$-module $E=M/\M M$. Then $M=\sum^s_{i=1}A\xi_i+\M M$,
and it follows from Corollary 4.3 that $M=\sum^s_{i=1}A\xi_i$.
Furthermore, the minimality of $\Omega$ follows from the gr-linear
independence of $\OV{\Omega}$ over $D$.\par\QED}\v5

Let $F=\oplus_{i\in J}Ae_i$ be a gr-free $A$-module with the
homogeneous free $A$-basis $\{ e_i\}_{i\in J}$ (see Section 1). Then
$\{ e_i\}_{i\in J}$ is certainly a minimal homogeneous generating
set of the $\Gamma$-graded $A$-module $F$. So, comparing with the
gr-rank of graded $D$-modules defined after Proposition 4.1, it
follows from the discussion of Section 3 that we may also define the
{\it gr-rank} of $F$ as the cardinality of $\{ e_i\}_{i\in J}$.
{\parindent=0pt\v5

{\bf Remark} At this stage, one should be aware of a delicate point,
namely we do not know generally whether an infinite minimal
homogeneous generating set of a gr-free $A$-module $F$ is a free
$A$-basis or not. But for a finitely generated gr-free $A$-module
$F$, this is not a problem (see Theorem 4.6 below).\v5

{\bf 4.5. Lemma} Let $M$ be a finitely generated $\Gamma$-graded
$A$-module with homogeneous generating set $\Omega =\{ \xi_1,\ldots
,\xi_s\}$. If $F=\oplus_{i=1}^sAe_i$ is the $\Gamma$-graded free
$A$-module of gr-rank $s$ with $d(e_i)=d(\xi_i)$, $\varphi$: $F\r M$
is the canonical $\Gamma$-graded epimorphism defined by $\varphi
(e_i)=\xi_i$, and $K=$ Ker$\varphi$, then $K\subseteq\M F$ if and
only if $\Omega$ is a minimal homogeneous generating set of
$M$.\vskip 6pt

{\bf Proof} Suppose that $\Omega$ is a minimal homogeneous
generating set of $M$. Note that $K$ is a $\Gamma$-graded submodule
of $F$. Let $f=\sum_{i=1}^sa_ie_i$ be a homogeneous element of
degree $\gamma$ in $K$. Then, by the construction of $F$, all the
$a_i$'s are homogeneous elements of $A$ and $d(a_1e_1)=\cdots
=d(a_se_s)=\gamma$. If $f\not\in\M F$, then since $A$ is
$\Gamma$-graded local, one of the coefficients, say $a_1$, must be
invertible. But $\varphi (f)=\sum^s_{i=1}a_i\xi_i=0$, so that
$\xi_1=-a_1^{-1}(\sum^s_{i=2}a_i\xi_i)$, contradicting the
minimality of $\Omega$.}\par

Conversely, suppose $K\subseteq \M F$. If $\Omega_1=\{
\xi_{i_1},\ldots ,\xi_{i_t}\}\subsetneq\Omega$ is such that
$M=\sum_{j=1}^tA\xi_{i_j}$ but some $\xi_q\not\in\Omega_1$, then by
Lemma 1.1, $\xi_q$ has a representation
$\xi_q=\sum_{j=1}^ta_j\xi_{i_j}$, in which the $a_j$'s are
homogeneous elements of $A$ such that $d(a_j\xi_{i_j})=d(\xi_q)$.
Thus, $e_q-\sum^t_{j=1}a_je_{i_j}\in$ Ker$\varphi =K\subseteq\M F$
and hence $e_q-\sum^t_{j=1}a_je_{i_j}=\sum^s_{i=1}b_ie_i$ with each
$b_i\in\M$ a homogeneous element. Comparing the coefficients of
$e_i$'s, it follows that $1\in\M$, a contradiction. This shows the
minimality of $\Omega$. \QED{\parindent=0pt \v5

{\bf 4.6. Theorem} Let $P$ be a finitely generated gr-projective
$A$-module (see Proposition 1.2). Then $P$ is a free $A$-module.
Indeed, any minimal homogeneous generating set $\Omega
=\{\xi_1,\ldots ,\xi_s\}$ of $P$ forms a free $A$-basis of
$P$.\vskip 6pt

{\bf Proof} Starting with a minimal homogeneous generating set
$\Omega =\{\xi_1,\ldots ,\xi_s\}$ of $P$, we have the gr-free
$A$-module $F=\oplus_{i=1}^sAe_i$ of gr-rank $s$ with
$d(e_i)=d(\xi_i)$, $1\le i\le s$, and the canonical $\Gamma$-graded
epimorphism $F\mapright{\varphi}{} P$  defined by $\varphi
(e_i)=\xi_i$. Put $K=$ Ker$\varphi$. Since $P$ is a projective
$A$-module, the exact sequence of $\Gamma$-graded modules and
$\Gamma$-graded $A$-homomorphisms
$$0~\mapright{}{}~K~\mapright{}{}~F~\mapright{\varphi}{}~P~\mapright{}{}~0$$
splits. So, $F=K\oplus P_0$ with $P_0\stackrel{\varphi}{\cong} P$.
By the minimality of $\Omega$ and Lemma 4.5, we have $K\subseteq \M
F$. Thus, $K\subseteq \M F=\M K\oplus\M P_0$ gives rise to $K=\M K$.
As $F$ is finitely generated, it follows from $K\cong F/P_0$ that
$K$ is finitely generated. Hence $K=0$ by Lemma 4.2. This shows that
$\varphi$ is an isomorphism, and consequently, $P$ is a free
$A$-module with $\Omega$ a free $A$-basis. \QED\v5

{\bf 4.7. Corollary} Suppose that $A$ is left gr-Noetherian and that
every finitely generated $\Gamma$-graded $A$-module has finite
projective dimension (i.e. $A$ is graded left regular). Then every
finitely generated $\Gamma$-graded $A$-module has a finite free
resolution, that is, there is an exact sequence of finite length
$$0~\mapright{}{}~F_n~\mapright{}{}~F_{n-1}~\mapright{}{}~\cdots~\mapright{}{}~F_0~\mapright{}{}~M~\mapright{}{}~0$$
in which each $F_i$ is a gr-free $A$-module of finite gr-rank.}\v5

\section*{5. The Existence of gr-Projective Covers and Minimal gr-Free Resolutions over
$A$} Let $\Gamma$ be a cancelation monoid with the neutral element
$e$, and let $A=\oplus_{\gamma\in\Gamma}A_{\gamma}$ be a
$\Gamma$-graded local ring in the sense of Definition 2.6.  As
before, we write $\M$ for the maximal graded ideal of $A$, and write
${\bf m}$ for the maximal ideal of the (classical) local ring $A_e$,
where ${\bf m}=\M\cap A_e$ by Proposition 2.4. Our aim in this
section is to establish, up to a graded isomorphism of chain
complex,{\parindent=.75truecm\par

\item{(1)} the existence of a unique gr-projective cover for every
finitely generated $\Gamma$-graded $A$-module, and \par

\item{(2)} the existence of a unique  minimal gr-free resolution for every finitely generated
$\Gamma$-graded $A$-module, in the case that $A$ is left
gr-Noetherian (i.e. each graded left ideal of $A$ is finitely
generated).\parindent=.5truecm\v5

Let $M$ be a $\Gamma$-graded module over the $\Gamma$-graded local
ring $A$. As in the ungraded case (e.g. [AF], P.72, P.199), we say
that a  graded submodule $N$ of $M$ is {\it gr-superfluous} if,
whenever $H$ is a graded submodule of $M$ with $H+N=M$, then $H=M$;
and we say that a $\Gamma$-graded epimorphism $\varepsilon$: $P\r
M$, where $P$ is a $\Gamma$-graded projective $A$-module, is a {\it
gr-projective cover} of $M$ if Ker$\varepsilon$ is
gr-superfluous.\par

Before proceeding to discuss the gr-projective cover, let us bear in
mind the fact that every gr-projective module is, as an ungraded
module, projective in the usual sense (Proposition 1.2). }
{\parindent=0pt\v5

{\bf 5.1. Proposition} The following two statements hold for a
$\Gamma$-graded $A$-module $M$. \par

(i) Suppose that $M$ has a gr-projective cover $\varepsilon$: $P\r
M$. Then for any $\Gamma$-graded epimorphism $\psi$: $Q\r M$, where
$Q$ is a $\Gamma$-graded projective $A$-module, there is a
$\Gamma$-graded epimorphism $\varphi$: $Q\r P$ such that the diagram
$$\begin{array}{ccccccc} &&&&Q&&\\
&&&\scriptstyle{\varphi}\swarrow&~\downarrow\scriptstyle{\psi}&&\\
&&P&\mapright{}{\varepsilon}&M&\rightarrow&0\\
&\swarrow&&&\downarrow{}&&\\
0&&&&0&&\end{array}$$ is commutative. Moreover,
Ker$\varphi\subseteq$ Ker$\psi$, and $Q=P_0\oplus$ Ker$\varphi$
with $P\cong P_0$.\par

(ii) If $f$: $M_1\r M_2$ is a graded isomorphism of $\Gamma$-graded
$A$-modules and if $\varepsilon_1$: $P_1\r M_1$ and $\varepsilon_2$:
$P_2\r M_2$ are gr-projective covers, then there is a graded
isomorphism of $\Gamma$-graded short exact sequences
$$\begin{array}{ccccccccc} 0&\r&\hbox{Ker}\varepsilon_1&\mapright{}{}&P_1&\mapright{\varepsilon_1}{}&M_1&\r&0\\
&&\mapdown{h}\scriptstyle{\cong}&&\mapdown{\varphi}\scriptstyle{\cong}&&\mapdown{f}\scriptstyle{\cong}&&\\
0&\r&\hbox{Ker}\varepsilon_2&\mapright{}{}&P_2&\mapright{\varepsilon_2}{}&M_2&\r&0\end{array}$$
Hence, any two gr-projective covers of $M$ are isomorphic. \vskip
6pt

{\bf Proof} (i) If $\psi$: $Q\r M$ is any $\Gamma$-graded
epimorphism, where $Q$ is a $\Gamma$-graded projective $A$-module,
then the projectiveness of $P$ and $Q$ yield two graded
homomorphisms $\varphi$ and $\psi$ such that the diagram
$$\begin{array}{ccccccc} &&&&Q&&\\
&&&\scriptstyle{\varphi}\swarrow\nearrow\scriptstyle{\varphi '}&~\downarrow\scriptstyle{\psi}&&\\
&&P&\mapright{}{\varepsilon}&M&\rightarrow&0\\
&&&&\downarrow{}&&\\
&&&&0&&\end{array}$$ is commutative in an obvious way, that is,
$\psi\varphi '=\varepsilon$, $\varepsilon\varphi =\psi$. It follows
that $\varepsilon (\varphi\varphi ')=\varepsilon$. Thus, for any
$x\in P$ we have $\varepsilon ((\varphi\varphi ')(x))=\varepsilon
(x)$, implying $x-\varphi (\varphi '(x))\in$ Ker$\varepsilon$. This
turns out that $P=\hbox{Im}\varphi +\hbox{Ker}\varepsilon .$ Since
Ker$\varepsilon$ is gr-superfluous,  $P=$ Im$\varphi$ and hence
$\varphi$ is epic. Concerning the remained assertions,
Ker$\varphi\subseteq$ Ker$\psi$ follows from $\varepsilon\varphi
=\psi$, and the direct decomposition $Q=P_0\oplus$ Ker$\varphi$ may
be obtained from the fact that the epimorphism $\varphi$ splits.
\par

(ii) By (i), there is a graded epimorphism $\varphi$: $P_1\r P_2$
such that $f\varepsilon_1=\varepsilon_2\varphi$. Since $f$ is an
isomorphism, it follows that $\hbox{Ker}\varphi =
\hbox{Ker}\varepsilon_1.$ On the other hand, since $\varphi$ splits
we have $P_1=$ Im$\psi\oplus$ Ker$\varphi =$ Im$\psi\oplus$
Ker$\varepsilon_1$, where $\psi$: $P_2\r P_1$ is the graded
homomorphism such that $\varphi\psi =1_{P_2}$. As Ker$\varepsilon_1$
is a gr-superfluous submodule of $P_1$, we have Im$\psi =P_1$. Thus,
$\psi$ is an isomorphism and so is $\varphi$. Finally, the desired
isomorphism of short exact sequences can be completed in an obvious
way.\QED }\v5

Based on Corollary 4.3, Lemma 4.5 and Proposition 5.1(ii), we are
able to derive the existence of a unique gr-projective cover for
every finitely generated graded $A$-module.{\parindent=0pt\v5

{\bf 5.2. Theorem} Every finitely generated $\Gamma$-graded
$A$-module $M$ has the following two properties.\par

(i) The graded submodule $\M M$ of $M$ is gr-superfluous.\par

(ii) Up to a graded isomorphism, $M$ has a unique gr-projective
cover.\par\QED }\v5

Let $M$ be a finitely generated $\Gamma$-graded $A$-module. We say
that $M$ has a {\it minimal gr-free resolution} if there is an exact
sequence of $\Gamma$-graded modules and graded homomorphisms
$$\cdots~\mapright{}{}~F_n~\mapright{\varphi_n}{}~F_{n-1}~\mapright{\varphi_{n-1}}{}~\cdots~
\mapright{\varphi_2}{}~F_1~\mapright{\varphi_1}{}~F_0~\mapright{\varphi_0}{}~M~\mapright{}{}~0$$
satisfying\par

(a) each $F_n$, $n\ge 0$, is a gr-free $A$-module of finite gr-rank,
and\par

(b) Im$\varphi_n\subseteq\M F_{n-1}$, $n\ge 1$.{\parindent=0pt\v5

{\bf 5.3. Theorem} Suppose that $A$ is left gr-Noetherian (i.e. each
graded left ideal is finitely generated). Then, up to a graded
isomorphism, every finitely generated $\Gamma$-graded $A$-module $M$
has a unique  minimal gr-free resolution.\vskip 6pt

{\bf Proof} Since $A$ is left gr-Noetherian by the assumption, a
minimal gr-free resolution of $M$ can be constructed by using Lemma
4.5 and Theorem 4.6; and the uniqueness of such a resolution is
guaranteed by Theorem 3.4 and Theorem 5.2.\QED }\v5

Let $\Gamma$ be a totally ordered cancelation monoid with the total
ordering $\prec$, such that the neutral element $e$ of $\Gamma$ is
the smallest element of $\Gamma$. If $M$ is a finitely generated
$\Gamma$-graded $A$-module over a $\Gamma$-graded local ring $A$,
then by Theorem 3.4(ii), any two minimal homogeneous generating sets
of $M$ contain the same number of homogeneous elements of degree
$\gamma$ for all $\gamma\in\Gamma$. Thus, if furthermore $A$ is left
gr-Noetherian and if
$${\cal
F}_{\bullet}\quad\quad\quad\cdots~\mapright{}{}~F_n~\mapright{\varphi_n}{}~F_{n-1}~\mapright{\varphi_{n-1}}{}~\cdots~
\mapright{\varphi_2}{}~F_1~\mapright{\varphi_1}{}~F_0~\mapright{\varphi_0}{}~M~\mapright{}{}~0$$
is a minimal gr-free resolution of $M$ in which each
$F_i=\oplus_{j=1}^{s_i}Ae_{ij}$ is a gr-free $A$-module with the
homogeneous free $A$-basis $\{ e_{ij}\}_{j=1}^{s_i}$, then, as with
a commutative polynomial algebra over a field $K$ or as with a
connected $\mathbb{N}$-graded algebra over a field $K$, we can
define, for each $\gamma\in\Gamma$, the {\it $\Gamma$-graded Betti
number} of $M$, denoted $\beta_{i,\gamma}(M)$,  as
$$\beta_{i,\gamma}(M)=\#\left\{\left. e_{it}\in\{e_{ij}\}_{j=1}^{s_i}~\right |~
d(e_{ij})=\gamma\right\} ;$$ and we can then define the {\it $i$-th
total $\Gamma$-graded Betti number} as
$\beta_i(M)=\sum_{\gamma\in\Gamma}\beta_{i,\gamma}(M)$. At this
stage, we will not further explore to what extent $A$ can be studied
by virtue of the $\Gamma$-graded Betti numbers, because such topic
has been beyond the scope of this paper. \v5

\section*{6. Determining Homological Dimensions via $A/\M$}

Let $A=\oplus_{\gamma\in\Gamma}A_{\gamma}$ be a $\Gamma$-graded
local ring (in the sense of Definition 2.6) with the graded maximal
ideal $\M$, where $\Gamma$ is a cancelation monoid with the neutral
element $e$. Based on the results obtained in previous sections, in
this section we demonstrate, by means of $A/\M$, how to determine
the projective dimension of a finitely generated $\Gamma$-graded
$A$-module $M$, and how to determine the (graded) left global 
homological dimension of $A$. \v5

We first fix some notations. Write $D$ for the $\Gamma$-graded
division ring $A/\M$, i.e. $D=A/\M$. Then $D_e\cong A_e/{\bf m}$ is
a division ring in the classical sense, where ${\bf m}=\M\cap A_e$
is the maximal ideal of the (classical) local ring $A_e$. For an
$A$-module $M$, we write p.dim$_AM$ for the projective dimension of
$M$ over $A$ in the usual sense. We also write gr.l.gl.dim$A$ for
the {\it graded left global dimension} of $A$, which is, in view of
Proposition 1.2, defined as
$$\hbox{gr.l.gl.dim}A=\sup\{ \hbox{p.dim}_AM~|~M~\hbox{any left}~\Gamma\hbox{-graded}~A\hbox{-module}\} .$$
Indeed, the definition gr.l.gl.dim given above needs only to
consider {\it all finitely generated} (or {\it all cyclic})
$\Gamma$-graded $A$-modules (see [LVO1], or see [LVO2],  where in
the proof of Proposition 2.5.3 it is enough to replace gr-free
modules over a group $G$-graded ring by gr-free modules over a 
monoid $\Gamma$-graded ring).  {\parindent=0pt\v5

{\bf 6.1. Lemma} Suppose that $A$ is left gr-Noetherian, i.e.,  each
graded left ideal of $A$ is finitely generated. Let $N$ be a nonzero
and finitely generated $\Gamma$-graded $A$-module. If
Tor$^A_1(D,N)=\{0\}$ then $N$ is a free $A$-module.\vskip 6pt

{\bf Proof} Since $N\ne \{ 0\}$ and is finitely generated, $E=N/\M
N\ne \{ 0\}$ by Lemma 4.2. It follows from Proposition 4.4 that $E$
is gr-free over $D=A/\M$, and if $\Omega =\{\xi_1,\ldots ,\xi_s\}$
is a minimal homogeneous generating set of $N$, then
$\OV{\Omega}=\{\OV{\xi_1},\ldots ,\OV{\xi_s}\}$ is a homogeneous
free $D$-basis for $E$, where each $\OV{\xi_i}=\xi_i+\M N$ is the
residue class of $\xi_i$ in $E$. Let $F=\oplus_{i=1}^sAe_i$ be the
gr-free $A$-module of gr-rank $s$ with the homogeneous free
$A$-basis $\{ e_i~|~d(e_i)=d(\xi_i),~1\le i\le s\}$, and consider
the short exact sequence
$$0~\mapright{}{}~K~\mapright{}{}~F~\mapright{\varphi}{}~N~\mapright{}{}~0$$
in which $\varphi (e_i)=\xi_i$, $1\le i\le s$, and $K=$
Ker$\varphi$. If Tor$^A_1(D,N)=\{ 0\}$, then it follows from
standard homological algebra (e.g. [Rot], P.188) that the sequence
$$0~\mapright{}{}~D\otimes_AK~\mapright{}{}~D\otimes_AF~\mapright{1\otimes\varphi}{}~D\otimes_AN~\mapright{}{}~0$$
is exact. But by the construction of $F$ and $\varphi$ we know that
$1\otimes\varphi$ is an isomorphism of gr-free $D$-modules. Hence
$D\otimes_AK=\{ 0\}$. Thus, the canonical isomorphism $K/\M K\cong
D\otimes_AK$ turns out that $K/\M K=\{ 0\}$. Note that $A$ is left
gr-Noetherian by our assumption. So, $K$ is finitely generated, and
by Lemma 4.2, $K=\{ 0\}$. This shows that $\varphi$ is an
isomorphism, i.e. $N$ is free, as desired.\QED \v5

{\bf 6.2. Theorem} Suppose that $A$ is left gr-Noetherian. The
following statements hold.\par

(i) Let $M$ be a finitely generated $\Gamma$-graded $A$-module.
Then, p.dim$_AM\le n<\infty$ if and only if Tor$^A_{n+1}(D, M)=0$
for some integer $n\ge 0$, where $D=A/\M$.\par

(ii) Let $M$ be a finitely generated $\Gamma$-graded $A$-module.
Then p.dim$_AM=$ the length of a minimal gr-free resolution of $M$
over $A$.\par

(iii) gr.l.gl.dim$A=$ p.dim$_AD$. \par

(iv) gr.l.gl.dim$A\le $ inj.dim$_AD$, where inj.dim$_AD$ stands for
the (classical) injective dimension of $D$ over $A$.\vskip 6pt

{\bf Proof} (i) If p.dim$_AM\le n<\infty$, then clearly
Tor$^A_{n+1}(D,M)=\{ 0\}$. Conversely, suppose that
Tor$^A_{n+1}(D,M)=\{ 0\}$ for some integer $n\ge 0$. Construct an
exact sequence of $\Gamma$-graded $A$-modules and graded
$A$-homomorphisms
$$0~\mapright{}{}~K_n~\mapright{}{}~F_{n-1}~\mapright{\varphi_{n-1}}{}~\cdots~\mapright{}{}~F_0~\mapright{\varphi_0}{}~M~\mapright{}{}~0$$
in which each $F_i$ is a gr-free $A$-module of finite gr-rank, and
$K_n=$ Ker$\varphi_{n-1}$ is finitely generated.  By standard
homological algebra (e.g. [Rot], Corollary 6.13), we have
$$\hbox{Tor}^A_1(D,K_n)=\hbox{Tor}^A_{n+1}(D,M)=\{ 0\} .$$
It follows from Lemma 6.1 that $K_n$ is free. Hence p.dim$_AM\le n$.
\par

(ii) Let
$${\cal F}_{\bullet}\quad\quad\quad\cdots~\mapright{}{}~F_n~\mapright{\varphi_n}{}~F_{n-1}~\mapright{\varphi_{n-1}}{}~\cdots~
\mapright{\varphi_2}{}~F_1~\mapright{\varphi_1}{}~F_0~\mapright{\varphi_0}{}~M~\mapright{}{}~0$$
be a minimal gr-free resolution of $M$ as constructed in Section 5.
Then since Im$\varphi_n\subseteq\M F_{n-1}$ for $n\ge 1$, the
derived sequences $D\otimes_A{\cal F}_{\bullet}$ and Hom$_A({\cal
F}_{\bullet},D)$ give rise to $$\hbox{Tor}^A_n(D,M)\cong
D^{i_n}\cong\hbox{Ext}_A^n(M,D),\quad\hbox{gr-rank}F_n=i_n,~n\ge
1.$$ Furthermore, by the well-known fact that Tor$^A_n(D,M)$ and
Ext$_A^n(M,D)$ are invariant with respect to projective resolution
we conclude that p.dim$_AM=$ the length of ${\cal F}_{\bullet}$.\par

(iii) Since $D=A/\M$ is a $\Gamma$-graded cyclic $A$-module, we have
p.dim$_AD\le$ gr.l.gl.dim$A$. Conversely, if p.dim$_AD=\infty$, then
gr.l.gl.dim$A\le $ p.dim$_AD$. In the case that p.dim$_AD=n<\infty$,
we have Tor$^A_{n+1}(D,M)=\{ 0\}$ for every finitely generated
$\Gamma$-graded $A$-module $M$. It follows from (i) that
gr.l.gl.dim$A\le n$. Thus, we have reached the equality
gr.l.gl.dim$A=$ p.dim$_AD$.\par

(iv)  If inj.dim$_AD=\infty$, then gr.l.gl.dim$A\le $ inj.dim$_AD$.
Suppose inj.dim$_AD=n<\infty$. We first prove the
following{\parindent=1.87truecm \par

\item{\bf Claim 1} If $N$ is a nonzero finitely generated
$\Gamma$-graded $A$-module with the property that $\M N=\{ 0\}$,
then inj.dim$_AN\le n$. \parindent=.5truecm

Let $\Omega =\{ \xi_1,...,\xi_s\}$ be a homogeneous generating set
of $N$. We prove Claim 1 by induction on $s$. If $s=1$ and
$d(\xi_1)=\gamma^*$, then the familiar $A$-module epimorphism
$A~\mapright{\varphi}{}~N=A\xi_1$ defined by $\varphi (a)=a\xi_1$
has the property that $\varphi (A_{\gamma})\subseteq
N_{\gamma\gamma^*}$ for all $\gamma\in\Gamma$. By means of the
cancelation law on $\Gamma$ it is easy to check that Ker$\varphi$ is
a graded left ideal of $A$. If $\M N=\{ 0\}$, then $\M\subseteq$
Ker$\varphi$ and consequently Ker$\varphi =\M$ because $\M$ is also
a maximal graded left ideal of $A$ by Section 2. Thereby $N\cong
A/\M=D$, and it follows that inj.dim$_AN=$ inj.dim$_AD=n$. If $s>1$,
then by the induction hypothesis we have inj.dim$_AN_1\le n$, where
$N_1=\sum^s_{i=2}A\xi_i$. Furthermore, the exact sequence
$$0~\mapright{}{}~N_1~\mapright{}{}~N~\mapright{}{}~N/N_1~\mapright{}{}~0$$
entails that inj.dim$_AN\le n$. So  Claim 1 is proved.}} \par

We complete the proof by showing the following{\parindent=0pt\par

{\bf Claim 2} If $M$ is a finitely generated $\Gamma$-graded
$A$-module, then p.dim$_AM\le n$.}\par

Construct an exact sequence of $\Gamma$-graded $A$-modules and
graded $A$-homomorphisms
$$0~\mapright{}{}~K~\mapright{}{}~F_{n-1}~\mapright{\varphi_{n-1}}{}~\cdots~\mapright{}{}~F_0~
\mapright{\varphi_0}{}~M~\mapright{}{}~0\eqno{(1)}$$ in which each
$F_i$ is a gr-free $A$-module of finite gr-rank, and $K=$
Ker$\varphi_{n-1}$ is finitely generated, and consider the exact
sequence
$$0~\mapright{}{}~L~\mapright{\ell}{}~F~\mapright{\psi}{}~K~\mapright{}{}~0\eqno{(2)}$$
in which $F$ is a gr-free $A$-module of finite gr-rank, $\psi$ is a
graded epimorphism, $L=$ Ker$\psi$ and $\ell$ is the inclusion map.
If $L=\{ 0\}$, then $K$ is free and (1) gives rise to p.dim$_AM\le
n$. If $L\ne \{ 0\}$, then since $A$ is left gr-Noetherian and hence
the graded submodule $L$ of $F$ is finitely generated, Lemma 4.2
ensures $N=L/\M L\ne \{ 0\}$. Since $\M N=\{ 0\}$, we have
inj.dim$_AN\le n$ by Claim 1. So if we consider an injective
resolution
$$0~\mapright{}{}~N~\mapright{}{}~Q_0~\mapright{}{}~Q_1~\mapright{}{}~\cdots~\mapright{}{}~Q_n~\mapright{}{}~0\eqno{(3)}$$
of $N$, then standard homological algebra tells us that
Ext$(M,Q_n)=\{ 0\}$ implies Ext$(K,N)=\{ 0\}$. It follows that the
exact sequence (2) yields the exact sequence
$$0~\mapright{}{}~\hbox{Hom}_A(K,N)~\mapright{}{}~\hbox{Hom}_A(F,N)~\mapright{}{}~\hbox{Hom}_A(L,N)~
\mapright{}{}~0\eqno{(4)}$$ Let $\pi$: $L\r N=L/\M L$ be the
canonical graded $A$-epimorphism. Then by the exactness of (4),
there is an $\alpha\in$ Hom$_A(F,N)$ mapping to $\pi$, namely the
diagram
$$\begin{array}{ccc} L&\mapright{\ell}{}&F\\
\mapdown{\pi}&\swarrow\scriptstyle{\alpha}&\\
N&&\\
\downarrow&&\\
0&&\end{array}\eqno{(5)}$$ is commutative. Note that $\ell$ and
$\pi$ are graded $A$-homomorphism. By using the cancelation law on
$\Gamma$, one may check that ([NVO], Lemma I.2.1) works with respect
to the diagram (5), i.e., there exists a graded $A$-homomorphism
$F~\mapright{\beta}{}~N$ such that $\beta\ell =\pi$ (see also the
proof of Proposition 1.2 in Section 1). Thus, by the gr-freeness of
$F$, there is a graded $A$-homomorphism $F~\mapright{\vartheta}{}~L$
such that $\pi\vartheta =\beta$. It follows from $\pi\vartheta\ell
=\beta\ell =\pi$ that $\pi (\vartheta\ell -1_L)=0$. Put $\rho
=\vartheta\ell -1_L$. Then $\pi\rho =0$ implies Im$\rho\subseteq$
Ker$\pi=\M L$. Moreover, since $x=(1_L+\rho )(x)-\rho (x)$ for all
$x\in L$, we have $L=$ Im$(1_L+\rho)+\M L$. So, by Lemma 4.2, $L=$
Im$(1_L+\rho )$. This shows that $1_L+\rho =\vartheta\ell$ is
surjective. Noticing that $L$ is a left gr-Noetherian $A$-module, we
may show, actually as in the ungraded case, that the graded
$A$-homomorphism $L~\mapright{\vartheta\ell}{}~L$ is an isomorphism.
Hence, the sequence (2) splits, i.e., there is an $A$-homomorphism
$F~\mapright{\tau}{}~L$ such that $\tau\ell =1_L$. This yields the
isomorphism $F\cong L\oplus K$, showing that $K$ is a projective
$A$-module. Turning back to the exact sequence (1), we conclude that
p.dim$_AM\le n$, as claimed.\QED\v5

The next theorem generalizes ([Li1], Theorem 3.1(2)) which deals 
with an {\it $\mathbb{N}$-graded} left Noetherian ring with $A_0$ a 
local ring in the classical sense.{\parindent=0pt\v5

{\bf 6.3. Theorem} Let $A$ be as in Theorem 6.2, and $D=A/\M$.
Assume further that $\Gamma$ is well-ordered by a well-ordering
$\prec$ such that the neutral element $e$ is the smallest element in
$\Gamma$. Then
$$\hbox{l.gl.dim}A=~\hbox{p.dim}_AD=~\hbox{inj.dim}_AD,$$ where
l.gl.dim$A$ is referred to as the left global dimension of
$A$.\vskip 6pt\def\F{{\cal F}}

{\bf Proof} Consider the $\Gamma$-grading filtration $\F A=\{
\F_{\gamma}A\}_{\gamma\in\Gamma}$ of
$A=\oplus_{\gamma\in\Gamma}A_{\gamma}$ defined by letting
$\F_{\gamma}A=\oplus_{\gamma '\preceq\gamma}A_{\gamma '}$,
$\gamma\in\Gamma$. Then $A$ is turned into a $\Gamma$-filtered ring
in the sense that $1\in \F_eA$,
$A=\cup_{\gamma\in\Gamma}\F_{\gamma}A$ and
$\F_{\gamma_1}A\F_{\gamma_2}A\subseteq \F_{\gamma_1\gamma_2}A$ for
all $\gamma_1,\gamma_2\in\Gamma$. Considering the associated
$\Gamma$-graded ring of $A$ with respect to $\F A$, which is by
definition the $\Gamma$-graded ring
$G(A)=\oplus_{\gamma\in\Gamma}G(A)_{\gamma}$ with
$G(A)_{\gamma}=\F_{\gamma}A/\F_{\prec\gamma}A$, where
$\F_{\prec\gamma}A=\cup_{\gamma '\prec\gamma}\F_{\gamma '}A$, it is
clear that $A\cong G(A)$ as $\Gamma$-graded rings. Hence, it follows
from ([Li2], Corollary 4.12) that l.gl.dim$A\le$ gr.lgl.dim$A$. But
it is clear that gr.l.gl.dim$A\le$ l.gl.dimn$A$. Therefore, Theorem
6.2 entails the desired equality
$\hbox{l.gl.dim}A=~\hbox{p.dim}_AD=~\hbox{inj.dim}_AD.$ }\v5

\centerline{References} \parindent=1.4truecm\par

\item{[AF]} F. W. Anderson and K. R. Fuller, {\it Rings and Categories of Modules}, Springer-Verlag, 1974.

\item{[BH]} W. Bruns and J. Herzog. {\it Cohen-Macaulay Rings}, Cambridge
Studies in Advanced Mathematics, Vol. 39. Cambridge University
Press, revised edition, 1998.

\item{[Eis]} D. Eisenbud, {\it Commutative Algebra with a View
Toward to Algebraic Geometry}, GTM 150, Springer, New York, 1995.

\item{[GW]} S. Goto and K. Watanabe, On graded rings I, {\it J. Math.
Soc. Japan}, 30(1978), 179--213.

\item{[Lam]} T. Y. Lam, {\it A First Course in Noncommutative Rings}, Graduate texts in Mathematics 131,
Springer-Verlag, 1991.

\item{[Lan]} S. Lang, {\it Algebra}, Third Edition, Springer-Verlag,
2002 (corrected in 2005).

\item{[Li1]} Huishi Li, Global dimension of graded local rings, {\it Comm. Alg.}, 24(7)(1996), 2399--2405.

\item{[Li2]} Huishi Li, $\Gamma$-Leading ~homogeneous algebras and
Gr\"obner bases, in: {\it Recent Developments in Algebra and Related
Areas} (F. Li and C. Dong eds.), Advanced Lectures in Mathematics,
Vol. 8, International Press \& Higher Education Press,
Boston-Beijing, 2009, 155 -- 200. {arXiv:math.RA/0609583}

\item{[LVO1]} Huishi Li and F. Van Oystaeyen, Global dimension and Auslander regularity of Rees rings, {\it Bull. Math. Soc. Belgique},
(serie A), Tom XLIII, 1991, 59--87.

\item{[LVO2]} Huishi Li and F. Van Oystaeyen, {\it Zariskian Filtrations}, $K$-monographs in Mathematics, Vol. 2,
Kluwer Academic Publishers, 1996.

\item{[MR]} J. C. McConnell and J. C. Robson, {\it Noncommutative Noetherian Rings}, J. Wiley, 1987.

\item{[NVO]} C.~N$\check{\rm a}$st$\check{\rm a}$sescu and F.~Van
Oystaeyen, {\it Graded Ring Theory}, Math. Library 28, North 
Holland, Amsterdam, 1982.

\item{[Rot]} J. J. Rotman, {\it An Introduction to Homological Algebra}, Academic Press, 1979.

\end{document}